\documentclass[10]{article}

\usepackage{amsmath,makeidx,amssymb,amscd,index,mathset,epsfig}
\usepackage[active]{srcltx}

\newcommand{\D}{\protect\displaystyle}
\newcommand{\T}{\protect\textstyle}
\newcommand{\eps}{\varepsilon}
\newcommand{\lbd}{\lambda}
\newcommand{\ipl}{\langle} 
\newcommand{\ipr}{\rangle} 
\newcommand{\eop}{ \raisebox{-.2ex}{\rule{2.5mm}{2.5mm}} }

\newtheorem{theorem}{Theorem}[section]

\newtheorem{lemma}[theorem]{Lemma}

\newtheorem{corollary}[theorem]{Corollary}

\title{\Large ON REGULARIZATION METHODS BASED ON DYNAMIC PROGRAMMING
TECHNIQUES}

\author{\small S. KINDERMANN$^{\mbox{\tiny 1}}$ \ and \
               A. LEIT\~AO$^{\mbox{\tiny2}}$}

\date{}

\begin{document}
\maketitle

\baselineskip=0.9\normalbaselineskip
\vspace{-3pt}

\begin{center}
{\footnotesize\em 
$^{\mbox{\tiny\rm 1}}$
Institute for Industrial Mathematics,
Johannes Kepler University, A-4040 Linz, Austria.
email: kindermann\symbol{'100}indmath.uni-linz.ac.at\\[3pt]
$^{\mbox{\tiny\rm 2}}$
Department of Mathematics, Federal University of St.\,Catarina,
88040-900 Florianopolis, Brazil.
email: aleitao\symbol{'100}mtm.ufsc.br}
\end{center}

\begin{abstract}
In this article we investigate the connection between {\em regularization
theory} for inverse problems and {\em dynamic programming} theory.
This is done by developing two new regularization methods, based on dynamic
programming techniques. The aim of these methods is to obtain stable
approximations to the solution of linear inverse ill-posed problems.
We follow two different approaches and derive a continuous and a discrete
regularization method.
Regularization properties for both methods are proved as well as rates of
convergence.
A numerical benchmark problem concerning integral operators with convolution
kernels is used to illustrate the theoretical results.
\end{abstract}



%\baselineskip=\normalbaselineskip

%---------------------------------------------------------------------------
% Section 1
%---------------------------------------------------------------------------
\section{Introduction}

Our main goal is to establish a connection between {\em regularization theory}
\cite{EHN96,Mor93} for inverse problems and {\em dynamic programming} theory
\cite{Bel53,Bel57,BDS62,BK65,Dre65} for optimal control problems of linear
quadratic type.
This is done by developing two new regularization methods, based on dynamic
programming results.
The first one is a continuous regularization method, derived from the
Hamilton-Jacobi Equation and the Pontryagin maximum principle. The second
is a discrete regularization method, derived from the Bellman optimality
principle.

In the sequel we describe the inverse problems we are concerned with.
Let $X$, $Y$ be Hilbert spaces. Consider the problem of finding $u \in X$
from the equation
\begin{equation} \label{eq:ip}
F u = y ,
\end{equation}
where $y \in Y$ represents the data and $F: X \to Y$ is a linear operator
modeling an ill-posed problem (e.g. a compact operator).
Since inverse of the operator $F$ is unbounded, the solution $u$ does not
depend in a stable way on the right hand side $y$ and regularization
techniques have to be used in order to obtain a stable solution.
Continuous and discrete regularization methods have been quite well
studied in the last two decades and one can find relevant information in
\cite{EHN96,EKN89,ES00,HNS95,Mor93,Tau94} and in the references therein.
The aim of these methods is to obtain stable approximations to the solution
of the inverse problem (\ref{eq:ip}).

Next we give a brief description of the optimal control problems (continuous
and discrete) that will serve as starting point for developing the
regularization methods in this article.
These problems are mainly characterized by possessing a linear dynamics and
a quadratic objective function.

Our first (continuous) approach is based on the the following constrained
optimization problem:
\begin{equation} \label{eq:ccp}
  \left\{ \begin{array}{l}
      {\rm Mimimize} \ J(x,w) := \int\limits_0^T
                     \ipl x(t), L x(t) \ipr + \ipl w(t), M w(t) \ipr \ dt \\
      {\rm s.t.} \\
      x' = Ax + Bw, \ t \ge 0 ,\ \ x(0) = x_0\, ,
   \end{array} \right.
\end{equation}
where $x(t) \in \mathbb R^n$ is the system trajectory, $w(t) \in
\mathbb R^m$, $t \ge 0$ is the control variable, $A, L \in \mathbb
R^{n,n}$, $B \in \mathbb R^{n,m}$, $M \in \mathbb R^{m,m}$ are given
matrices and $x_0 \in \mathbb R^n$ is the initial condition.
The goal of the control problem is to find a pair of functions $(x,w)$, 
minimizing the quadratic objective function $J$ and satisfying the 
constraint imposed by the linear dynamical system -- such pairs are 
called {\em admissible processes}. In this article we adapt a solution
technique for this problem (dynamic programming) in order to derive a
continuous regularization method for the inverse problem (\ref{eq:ip}).

% This is achieved as follows: 
% First we have to define an optimal control problem related to (\ref{eq:ip}),
% what is done by choosing an objective function related to the residual
% $\| F(u) - y \|$. Furthermore, we use as initial condition any 
% approximation $u_0 \in X$ for the least square solution $u^\dag$
% of (\ref{eq:ip}). The dynamics corresponds to the choice of a descent
% direction for the regularization method. This completes the definition
% of the optimal control problem.
% 
% The next step in the formulation of the regularization method consists 
% in using the spectral decomposition of the operator $F$, in order to 
% obtain, from the dynamic programming approach, an optimal process for 
% the control problem in (\ref{eq:ccp}). Finally, we prove that the 
% optimal trajectories $\bar u(t)$, $t \in [0,T]$, generate a family 
% of regularization operators $R_T := u(T)$ for problem (\ref{eq:ip}),
% in the sense of \cite{EHN96}.

Our second (discrete) approach, has as starting point the discrete optimal
control problem
\begin{equation} \label{eq:dcp}
  \left\{ \begin{array}{l}
      {\rm Mimimize} \ J(x,w) := \ipl x_N , S x_N \ipr +
      \sum\limits_{k=0}^{N-1} \ipl x_k , L x_k \ipr +
                                 \ipl w_k , M w_k \ipr \\
      {\rm s.t.} \\
      x_{k+1} = A x_k + B w_k, \ k=0,\dots,N-1,\ \ x_0 \in \mathbb R^n\, .
   \end{array} \right.
\end{equation}
The matrices $A$, $B$, $L$, $M$ have the same meaning as in problem
(\ref{eq:ccp}) and $S \in \mathbb R^{n,n}$ is positive definite. Notice
that the final time $T$ in (\ref{eq:ccp}) is substituted by the number
of discrete steps $N \in \mathbb N$ in (\ref{eq:dcp}).
Again, using the dynamic programming technique, we are able to derive
a discrete regularization method for the inverse problem (\ref{eq:ip}).
In this discrete framework, the dynamic programming approach consists
basically of the Bellman optimality principle and the dynamic programming
equation.

To the best of our knowledge, dynamic programming techniques have only
been applied to solve particular inverse problems so far.
In \cite{KS97} the inverse problem of identifying the initial condition
in a semilinear parabolic equation is considered. In \cite{KS99} the same
authors consider a parameter identification problem for identification of
systems of distributed parameters.
In the engineering literature, dynamic programming filter has been used as
a regularization tool for inverse problems \cite{TB97,TB83}.
In this article however, dynamic programming methods allow us to formulate
regularization methods in an abstract functional analytical framework for
general inverse problems.

The article is outlined as follows: In Section~\ref{sec:deriv} we 
derive both regularization methods (continuous and discrete).
In Section~\ref{sec:regul} we analyze regularization properties of the
proposed methods. Rates of convergence are derived under abstract source
conditions and an {\em a priori} parameter choice yielding optimal order
convergence rates is provided.
Furthermore, for the discrete regularization method, we characterize
the filter functions (for the regularization operator) in terms of 
Chebyshev polynomials.
In Section~\ref{sec:num} numerical realizations of our regularization
methods are presented. We use our methods to solve an integral equation
of the first kind and compare the obtained performances with the Landweber
iteration and with the CG-method.

%---------------------------------------------------------------------------
% Section 2
%---------------------------------------------------------------------------
\section{Derivation of the regularization methods} \label{sec:deriv}

%---------------------------------------------------------------------------
\subsection{A continuous approach}

We start this section defining an optimal control problem related
with the linear inverse problem (\ref{eq:ip}). Let
$u_0 \in X$ be any approximation for the minimum norm solution
$u^\dag \in X$ of (\ref{eq:ip}). We aim to find a function
$u: [0,T] \to X$ such that, $u(0) = u_0$ and
\begin{equation} \label{eq:ip-asymp}
\| F u(T) - y \| \ \approx \ \| F u^\dag - y \|
\end{equation}
(notice that, if the data are contaminated with noise, it may not belong
to the range of $F$).
In the control literature, the function $u$ is called {\em trajectory}
(or {\em state}) and its evolution is is described by a dynamical system.
For simplicity, we choose a linear evolution model, i.e.
$u'  =  A u(t) + B v(t)$, $t \ge 0$,
where $A, B : X \to X$ are linear operators and $v:
[0,T] \to X$ is the {\em control} of the system (compare
with the classical problem in (\ref{eq:ccp})).
Keeping in mind the desired property described in (\ref{eq:ip-asymp}),
we shall see that for the proof of the convergence and stability results
in Section~\ref{sec:regul} it is enough to consider a simpler dynamics,
which does not depend on the state $u$, but only on the control $v$.
This justifies the choice of the dynamics: $u' = v$, $t \ge 0$.
In this case, the control $v$ corresponds to a {\em velocity function}.

The next step is to choose the objective function for our control problem.
Recalling the formulation of the linear quadratic control problem in
(\ref{eq:ccp}) and also the goals described in (\ref{eq:ip-asymp}), the
objective function has to be related to the minimization of both the
residual norm and the velocity norm along the trajectories, i.e.
$$ J(u,v) := \T\frac{1}{2} \D\int_0^T \|F u(t) - y\|^2 + \| v(t) \|^2\ dt . $$
Putting all together we obtain the following abstract optimal control
problem in Hilbert spaces:
\begin{equation} \label{eq:ccp-hs}
  \left\{ \begin{array}{l}
      {\rm Mimimize} \ J(u,v) = \frac{1}{2} \displaystyle
                     \int_0^T \| F u(t) - y \|^2 + \| v(t) \|^2 \ dt \\
      {\rm s.t.} \\
      u' = v, \ t \ge 0 ,\ \ u(0) = u_0\, ,
   \end{array} \right.
\end{equation}
where the (fixed but arbitrary) final time $T > 0$ will play the role
of the regularization parameter. The functions $u,v : [0,T] \to X$
correspond respectively to the trajectory and the control of the system.
The pairs $(u,v)$ are called {\em processes}.

Next we define the residual function $\eps(t) := F u(t) - y$ associated
to a given trajectory $u$. Notice that this residual function evolves
according to the dynamics
$$ \eps' = F u'(t) = F v(t)\, ,\ t \ge 0\, . $$
With this notation, problem (\ref{eq:ccp-hs}) can be rewritten in the
following form
\begin{equation} \label{eq:ccpr-hs}
  \left\{ \begin{array}{l}
      {\rm Mimimize} \ J(\eps,v) = \frac{1}{2} \D
                     \int_0^T \| \eps(t) \|^2 + \| v(t) \|^2 \ dt \\
      {\rm s.t.} \\
      \eps' = F v, \ t \ge 0 ,\ \ \eps(0) = F u_0 - y\, .
   \end{array} \right.
\end{equation}

It is immediate to detect a parallel between solvability of the
optimal control problem (\ref{eq:ccp-hs}) and the auxiliary problem
(\ref{eq:ccpr-hs}), namely:
if $(\bar u, \bar v)$ is an optimal process for problem
(\ref{eq:ccp-hs}), then the process $(\bar \eps, \bar v)$, with
$\bar \eps := F \bar u - y$, will be an optimal process for problem
(\ref{eq:ccpr-hs}). Conversely, if $(\bar \eps, \bar v)$ is an optimal
process for problem (\ref{eq:ccpr-hs}), with $\eps(0) = F u_0 - y$,
for some $u_0 \in X$, then the corresponding process $(\bar u, \bar v)$
is an optimal process for problem (\ref{eq:ccp-hs}).

In the sequel, we derive the dynamic programming approach for the
optimal control problem in (\ref{eq:ccpr-hs}). We start by introducing
the first Hamilton function. This is the function $H: \mathbb R \times
Y^2 \times X \to \mathbb R$ given by
$$  H(t,\eps,\lbd,v) \ := \ \D\ipl \lbd , F v \ipr +
    \T\frac{1}{2}\D [ \ipl \eps , \eps \ipr + \ipl v , v \ipr ] \, . $$
Notice that the variable $\lbd$ plays the role of a Lagrange multiplier
in the above definition. According to the Pontryagin's maximum principle,
the Hamilton function furnishes a necessary condition of optimality for
problem (\ref{eq:ccpr-hs}). Furthermore, since this function (in this
particular case) is convex in the control variable, this optimality
condition also happens to be sufficient. Recalling the maximum principle,
along an optimal trajectory we must have
\begin{equation} \label{eq:opt-cond}
0 \ = \ \frac{\partial H}{\partial v}(t,\eps(t),\lbd(t),v(t))
      \ = \ F^* \lbd(t) + v(t) \, .
\end{equation}
This means that the optimal control $\bar v$ can be obtained directly
from the Lagrange multiplier $\lbd: [0,T] \to Y$, by the formula
$$  \bar v(t) = -F^* \lbd(t)\, ,\ \forall t\, . $$
Therefore, the key task is actually the evaluation of the Lagrange
multiplier. This leads us to the Hamilton-Jacobi equation. Substituting
the above expression for $\bar v$ in (\ref{eq:opt-cond}), we can define
the second Hamilton function $\mathcal H: \mathbb R \times Y^2
\to \mathbb R$
$$  \mathcal H(t,\eps,\lbd) \ := \ \min_{v \in X}
    \{ H(t,\eps,\lbd,v) \} \ = \ \T\frac{1}{2} \D\ipl \eps , \eps \ipr -
    \T\frac{1}{2} \D\ipl \lbd , F F^* \lbd \ipr \, . $$
Now, let $V: [0,T] \times Y \to \mathbb R$ be the value function
for problem (\ref{eq:ccpr-hs}), i.e.
\begin{eqnarray}
  V(t,\xi) \!\!\! & := & \!\!\!
                  \min\Big\{ \T\frac{1}{2} \int_t^T \| \eps(s) \|^2 +
                  \|v(s)\|^2 \, ds \ \Big| \ (\eps,v) \
                  {\rm admissible\ process} \nonumber \\
           \!\!\! &    & \!\!\!
                  {\rm \ \ \ \ \ \ \ \ for\ problem\
                  (\ref{eq:ccpr-hs})\ with\ initial\ condition} \
                  \eps(t) = \xi \Big\} \, . \label{def:cvf}
\end{eqnarray}
The interest in the value function follows
from the fact that this function is related to the Lagrange multiplier
$\lbd$ by the formula: $\lbd(t) = \partial V / \partial\eps (t,\bar \eps)$,
where $\bar \eps$ is an optimal trajectory.

From the control theory we know that the value function is a solution
of the Hamilton-Jacobi equation
\begin{equation} \label{eq:hjb}
 \frac{\partial V}{\partial t}(t,\eps) +
  \mathcal H(t,\eps,\frac{\partial V}{\partial \eps}(t,\eps)) \ = \ 0\, .
\end{equation}
Now, making the ansatz: $V(t,\eps) = \frac{1}{2} \ipl \eps , Q(t) \eps \ipr$,
with $Q: [0,T] \to L(Y,Y)$, we are able to rewrite (\ref{eq:hjb}) in the form
$$  \ipl \eps , Q'(t) \eps \ipr + \ipl \eps , \eps \ipr - 
    \ipl Q(t) \eps , F F^* Q(t) \eps \ipr \ = \ 0 \, . $$
Since this equation must hold for all $\eps \in X$, the function
$Q$ can be obtained by solving the Riccati equation
\begin{equation} \label{eq:riccati}
Q'(t) \ = \ -I + Q(t) F F^* Q(t) \, .
\end{equation}
Notice that the cost of all admissible processes for an initial
condition of the type $(T,\eps)$ is zero. Therefore we have to
consider the Riccati equation (\ref{eq:riccati}) with the final
condition
\begin{equation} \label{eq:riccati-ic}
Q(T) \ = \ 0\, .
\end{equation}

Once we have solved the initial value problem (\ref{eq:riccati}),
(\ref{eq:riccati-ic}), the Lagrange multiplier is given by $\lbd(t)
= Q(t) \bar\eps(t)$ and the optimal control is obtained by the
formula $\bar v(t) = - F^* Q(t) \bar\eps(t)$. Therefore, the optimal
trajectory of problem (\ref{eq:ccp-hs}) is defined via
\begin{equation} \label{eq:cont-reg}
 \bar u' = - F^* Q(t) [F \bar u(t) - y] \, ,\ \ \bar u(0) = u_0 \, .
\end{equation}

We use the optimal trajectory defined by the initial value problem 
(\ref{eq:cont-reg}) in order to define a family of reconstruction
operators $R_T: X \to X$, $T \in \mathbb R^+$,
\begin{equation} \label{eq:cont-Rt}
R_T(y) \ := \ \bar u(T) \ = \ u_0 -
              \int_0^T F^* Q(t) [F \bar u(t) - y] \ dt\, .
\end{equation}
We shall return to the operators $\{ R_T \}$ in Section~\ref{sec:regul}
and prove that the family of operators defined in (\ref{eq:cont-Rt})
is a regularization method for (\ref{eq:ip}) (see, e.g.,
\cite[Section~3.1]{EHN96}).

%---------------------------------------------------------------------------
\subsection{A discrete approach}\label{discapp}

In this section we use the optimal control problem (\ref{eq:dcp}) as
starting point to derive a discrete reconstruction method for the
inverse problem in (\ref{eq:ip}). Again, let $u_0 \in X$ be a given
approximation for the minimum norm solution $u^\dag \in X$
of (\ref{eq:ip}) and $N \in \mathbb N$. Analogously as we did in the
previous section, we aim to find a sequence $\{ u_k \}_{k=1}^N$ in
$X$, starting from $u_0 = u_0$, such that
\begin{equation} \label{eq:ip-asymp-d}
\| F u_N - y \| \ \approx \ \| F u^\dag - y \| \, .
\end{equation}
As in the previous section, we have now a discrete trajectory,
represented by the sequence $u_k$, which evolution is described by
the discrete dynamics
$$ u_{k+1} \ = \ A u_k \, + \, B v_k\, , \ k =0,1,\dots  $$
where the operators $A$ and $B$ are defined as before and $\{ v_k \}_{k=0}^
{N-1}$, is the control of the system (compare with (\ref{eq:dcp})).
As in the continuous case, we shall consider a simpler dynamics:
$u_{k+1} = u_k + v_k$, $k=0,1,\dots$ (i.e., $A = B = I$). To simplify
the notation, we represent the processes $(u_k, v_k)_{k=1}^N$ by $(u,v)$.

The objective function is chosen similarly as in the continuous case:
$$ J(u,v) := \T\frac{1}{2} \D \ipl F u_N - y ,\, S (F u_N - y) \ipr 
   + \T\frac{1}{2} \sum\limits_{k=0}^{N-1} \|F u_k-y\|^2 + \|v_k\|^2\, , $$
with some positive operator $S: Y \to Y$. Putting all
together we obtain the following abstract optimal control problem in
Hilbert spaces:
\begin{equation} \label{eq:dcp-hs}
  \left\{ \begin{array}{l}
      {\rm Mimimize} \ \ J(u,v) = \frac{1}{2}
                     \ipl F u_N - y ,\, S (F u_N - y) \ipr \\
      {\hskip3.7cm}     + \frac{1}{2} \sum_{k=0}^{N-1} \|F u_k-y\|^2 +
                     \|v_k\|^2 \\
      {\rm s.t.} \\
      u_{k+1} = u_k + v_k, \ k=0,1,\dots ,\ \ u_0 \in X
   \end{array} \right.
\end{equation}
where the (fixed but arbitrary) number of discrete steps $N \in
\mathbb N$ will play the role of the regularization parameter.

As in the continuous approach, we define the residual sequence $\eps_k :=
F u_k - y$, associated to a given trajectory $u$. Notice that
$$ \eps_{k+1} = F u_{k+1} - y = \eps_k + F v_k\, ,\ k = 0,1,\dots  $$
With this notation, problem (\ref{eq:dcp-hs}) can be rewritten in the form
\begin{equation} \label{eq:dcpr-hs}
  \left\{ \begin{array}{l}
      {\rm Mimimize} \ J(\eps,v) = \frac{1}{2} \ipl \eps_N, S\eps_N \ipr
      + \frac{1}{2} \sum\limits_{k=0}^{N-1} \| \eps_k \|^2 + \| v_k \|^2 \\
      {\rm s.t.} \\
      \eps_{k+1} = \eps_k + F v_k, \ k=0,1,\dots ,\ \ \eps_0 = F u_0 - y\, .
   \end{array} \right.
\end{equation}

Notice that if $(\bar u, \bar v)$ is an optimal process for problem 
(\ref{eq:dcp-hs}), then the process $(\bar \eps, \bar v)$, with 
$\bar \eps_k := F \bar u_k - y$, will be an optimal process for
problem (\ref{eq:dcpr-hs}) and vice versa, as one can easily check.

In the sequel, we derive the dynamic programming approach for the
optimal control problem in (\ref{eq:dcpr-hs}). We start by introducing
the value function (or Lyapunov function) $V: \mathbb R \times Y
\to \mathbb R$,
$$  V(k,\xi) \ := \ \min\{ J_k(\eps,v)\; |\ (\eps,v)
                     \in Z_k(\xi) \times X^{N-k} \} \, , $$
where
$$  J_k(\eps,v) \ := \ \T\frac{1}{2} \Big[ \ipl \eps_N , S \eps_N \ipr
    + \T\sum\limits_{j=k}^{N-1} \| \eps_j \|^2 + \| v_j \|^2 \Big]  $$
and
$$  Z_k(\xi) \ := \ \{ \eps \in Y^{N-k+1}\; |\ \eps_k = \xi ,\
                   \eps_{j+1} = \eps_j + F 
v_j,\ j=k,\dots,N-1 \} \, . $$
(Compare with the definition in (\ref{def:cvf})). The Bellman
principle for this discrete problem reads
\begin{equation} \label{eq:hjb-d}
V(k,\xi) \ = \ \min\{ V(k+1, \xi + Fv) + \T\frac{1}{2}
    ( \ipl \xi , \xi \ipr + \ipl v , v \ipr)\; |\ v \in X\} \, .
\end{equation}
The optimality equation (\ref{eq:hjb-d}) is the discrete counterpart of
the Hamilton-Jacobi equation (\ref{eq:hjb}). Notice that the value
function also satisfies the boundary condition: \ $V(N,\xi) \ = \ 
\frac{1}{2} \ipl \xi, S \xi \ipr$.

As in the continuous case, the optimality equation have to be solved
backwards in time ($k = N-1, \dots, 1$) recursively.

For $k=N-1$, we have
\begin{equation} \label{eq:ow_tm1}
  V(N-1,\xi) \ = \ \min\{ \T\frac{1}{2} ( \ipl \xi+Fv ,\, S(\xi+Fv) \ipr
           + \ipl \xi , \xi \ipr + \ipl v , v \ipr)\; |\ v \in X \} .
\end{equation}
A necessary and sufficient condition for $u_{N-1}$ to be a minimum of
(\ref{eq:ow_tm1}) is given by \ $v + F^* S(\xi + Fv) \ = \ 0$.
Solving this equation for $v$ we obtain
$$  \bar v_{N-1} \ := \ -(F^* S F + I)^{-1} F^* S \xi \, . $$
In order to obtain the optimal control recursively, we evaluate the matrices
\begin{equation}\label{eq:dcs}
 \begin{array}{l}
      S_N \ := \ S ; \\
      \mbox{for \ $k = N-1, \ldots, 0$ \ evaluate} \\
         \hskip1cm  R_k \ := \ (F^* S_{k+1} F + I)^{-1} F^* S_{k+1} \, ; \\
         \hskip1cm  S_k \ := \ (I-FR_k)^* S_{k+1} (I-FR_k) + R_k^* R_k + I\, ;
   \end{array} 
\end{equation}

\noindent
Once the matrices $R_k$ and $S_k$ are known, we obtain the optimal control
recursively, using the algorithm: \medskip
\begin{equation}\label{eq:dcu}
 \begin{array}{l}
      \eps_0 := F u_0 - y\, ; \\
      \mbox{for \ $k = 0, \ldots, N-1$, \ evaluate} \\
         \hskip1cm  \bar v_k \, := \, - R_k \bar \eps_k\, ; \\
         \hskip1cm  \bar u_{k+1} \ := \ u_k + \bar v_k \, ; \\
         \hskip1cm  \bar \eps_{k+1} \ := \ \bar \eps_k + F \bar v_k \, ; \\
  \end{array} 
\end{equation}

\medskip \noindent
to obtain the optimal control $\bar v = (\bar v_0, \ldots, \bar v_{N-1})$,
the optimal trajectory for problem (\ref{eq:dcpr-hs}) $\bar \eps =
(\bar \eps_0, \ldots, \bar \eps_N)$, and the optimal trajectory for
problem (\ref{eq:dcp-hs}) $\bar u = (\bar u_0, \ldots, \bar u_N)$.
Furthermore, the optimal cost is given by $V(0,\eps_0) = \frac{1}{2}
\ipl \eps_0 , S_0 \eps_0 \ipr$.

%---------------------------------------------------------------------------
% Section 3
%---------------------------------------------------------------------------
\section{Regularization properties} \label{sec:regul}

%---------------------------------------------------------------------------
\subsection{Regularization in the continuous case}

In this section we investigate the regularization properties
of the operator $R_T$ introduced in (\ref{eq:cont-Rt}).
Consider the Riccati equation (\ref{eq:riccati}) for the operator $Q$: 
We may express the operator $Q(t)$ via the spectral
family $\{F_\lambda\}$ of $FF^*$ (see e.g. \cite[Section~2.3]{EHN96}).
Hence, we make the ansatz
$$ Q(t) = \int q(t,\lambda) d F_\lambda \; . $$
Assuming that $q(t,\lambda)$ is $C^1$ we may find
from (\ref{eq:riccati}) together with the boundary condition at $t = T$
that
$$  \int \left( \T\frac{d}{d t} q(t,\lambda)  
    + 1 -  q(t,\lambda)^2 \lambda  \right) \ d F_\lambda  = 0, \ \ \ 
    q(T,\lambda) = 0 . $$
Hence, we obtain an ordinary differential equation  for $q$:
\begin{equation}\label{ode}
 \T\frac{d}{d t} q(t,\lambda) = -1 + \lambda q(t,\lambda)^2
\end{equation}
The solution to these equations is given by
\begin{equation}\label{qfunc}
 q(t,\lambda) = - \frac{1}{\sqrt{\lambda}} 
\tanh(\sqrt{\lambda}(t -T)) = 
 \frac{1}{\sqrt{\lambda}} 
\tanh(\sqrt{\lambda}(T-t)) .
\end{equation}
If $t < T$, then $Q(t)$ is nonsingular, since
$\lim_{x\to 0 } \frac{\tanh(x a)}{x}  = a$ and 
$\frac{\tanh( a x)}{x} $ is monotonically decreasing 
for $x >0$. Hence the spectrum of $Q(t)$ is
contained in the interval $[\frac{\tanh((T-t)\|F\|)}{\|F\|}, (T-t)]$.
Now consider the evolution equation (\ref{eq:cont-reg}):
The operator $Q(t)$ can be expressed as
$Q(t) = q(t,FF^*)$;
by usual spectral theoretic properties  (see, e.g., \cite[Page~44]{EHN96})
it holds that
$$ F^*q(t,F F^*)  =  q(t,F^* F)F^*. $$
Hence we obtain the problem
\begin{eqnarray}
 u'(t) &=&  -q(t,F^*F) \left( F^*F u(t) - F^*y \right)  \label{ip1}\\
 u(0) &=& u_0 \label{ip2}
\end{eqnarray}
%

% By linearity we may write  solution to the initial value problem
% (\ref{ip1},\ref{ip2}) 
% $u = u_h + u_i$,
% where $u_h$ solves (\ref{ip1}) with the initial condition
% $u_h(0) = 0$ and  
% $u_i$ solves
% \begin{eqnarray}\label{initial}
% u'(t) &=&  -q(t,F^*F) F^*F u(t)  \\
% u(0) & = & u_0 
% \end{eqnarray}

We may again use an ansatz via spectral calculus: if we set
$$ u(t) = \int g(t,\lambda) d E_\lambda F^* y$$
where $E_\lambda$ is the spectral family of $F^* F$,
we derive an ordinary differential equation for $g$.
%
%$$ \int_{\sigma} \frac{d}{dt} g(t,\lambda) 
%d E_\lambda F^*y = 
%\int_{\sigma} - q(t,\lambda) \lambda g(t,\lambda) 
%+ q(t,\lambda) d E_\lambda F^*y $$
%Hence
%$$ g'(t,\lambda) = - q(t,\lambda) \lambda g(t,\lambda) + q(t,\lambda)$$
%with $q$ as in (\ref{qfunc})
%
Similar as above, we can express the solution to (\ref{ip1},\ref{ip2}) 
in the form
%
%$g(0,\lambda) = 0$  is given by
%\begin{equation}
% g(t,\lambda) =  \frac{1 - \frac{\cosh(\sqrt{\lambda}(T-t))}{
%\cosh(\sqrt{\lambda} T )}}{\lambda}
%\end{equation}
%Similar we find for $u_i$,
%$$ u_i = \int_{\sigma} h(t,\lambda) d E_\lambda u_0$$
%with 
%$$ h(t,\lambda) = \frac{\cosh(\sqrt{\lambda}(T-t))}{\cosh(\sqrt{\lambda}T )}
%$$
%
\begin{equation}\label{tsol}
u(t) = 
\int \frac{1 - \frac{\cosh(\sqrt{\lambda}(T-t))}{
\cosh(\sqrt{\lambda} T)}}{\lambda} d E_\lambda F^*y 
+ \int
 \frac{\cosh(\sqrt{\lambda}(T-t))}{\cosh(\sqrt{\lambda} T )} d E_\lambda u_0.
\end{equation}
Setting $t=T$ we find an approximation of the solution 
\begin{equation}\label{solution}
u_T:=  u(T) = 
\int \frac{1 - \frac{1}{
\cosh(\sqrt{\lambda} T )}}{\lambda} d E_\lambda F^*y 
+ \int
 \frac{1}{\cosh(\sqrt{\lambda} T )} d E_\lambda u_0 .
\end{equation}
Note the similarity to Showalter`s methods \cite[Page~77]{EHN96},
where the term $\exp(\lambda T)$ instead of $\cosh( \sqrt{\lambda} T )$
appears.

\begin{theorem} \label{th:convcont}
The operator $R_T$ in (\ref{eq:cont-Rt}) is a regularization operator with
qualification $\mu_0 = \infty$ \cite[Page~76]{EHN96}, i.e. it satisfies

i) If the data are exact, $y = F u^\dagger$ and 
$u^\dagger$ satisfies a source condition for some
$\mu > 0$ 
\begin{equation}\label{sourcecond} 
 \exists\ \omega \in X : \quad  u^\dagger  = (F^*F)^\mu \omega ,
\end{equation}
we have the estimate
$$ \|u_T - u^\dagger\| \leq C_\mu T^{-2 \mu} $$

ii) If the data are contaminated with noise,
$ \| y -y_\delta\| \leq \delta $ and $y = F u^\dagger$ 
with $u^\dagger$ as in  (\ref{sourcecond}), then we have
$$ \|u_{T,\delta} - u^\dagger\| \leq
C_\mu T^{-2 \mu} + \delta T . $$
In particular, the a-priori parameter choice
$T \sim \delta^{\frac{-1}{2 \mu +1}} $ yields the optimal
order convergence rate
$$ \|u_{T,\delta} - u^\dagger\| \sim \delta^{\frac{2 }{2 \mu +1}} . $$
\end{theorem}

\noindent {\em Proof:}
For simplicity we set $u_0 = 0$, the generalization to 
the inhomogeneous case is obvious.
 (\ref{solution}) gives an expression of 
the regularization operator in terms of a 
filter function:
$$ R_T = \int f(T,\lambda) d E_{\lambda} F^* y $$
with 
$$ f(T,\lambda) = \lambda^{-1} \ \left( 1 - \frac{1}{\cosh(\sqrt{\lambda} T
)} \right) . $$
According to \cite[Theorem~4.1]{EHN96} we have to show that the filter
function $f(T,\lambda)$ satisfies the properties (regarding $1/T$ as 
regularization parameter).
\begin{enumerate}
\item for $T$ fixed, $f(T,.)$ is continuous;
\item there exists a constant $C$ such that
for all $\lambda >0$
$$|\lambda f(T,\lambda) | \leq C ; $$
\item $\hfil \lim\limits_{T \to \infty} f_T(\lambda) = \lambda^{-1}\, ,\ \ 
\forall \lambda \in (0, \| F^*F \|] . \hfil $
\end{enumerate}

\noindent
1. is clear since 
$ \lim_{\lambda \to 0} f(T,\lambda) = \frac{T^2}{2} $
the function can be extended continuously to $\lambda = 0$.

\noindent
2.  holds with $C = 1$ since
$ 0 \leq \frac{1}{\cosh(\sqrt{\lambda}(T))} \leq 1 $.

\noindent
3.  is obviously is the case since $\lim_{T\to\infty}\cosh(s) = \infty$.

We have to show that the qualification $\mu_0 = \infty$:
this needs an estimate $w_\mu(T)$ such that
$$ \lambda^\mu |(1- \lambda f(T,\lambda)|) \leq w_\mu(T). $$
It holds that
$$  \lambda^\mu |(1- \lambda f(T,\lambda)|)
= \frac{\lambda^{\mu}}{\cosh(\sqrt{\lambda} T)} \leq
2\frac{\lambda^{\mu}}{\exp(\sqrt{\lambda} T)} \leq
2 (2 \mu)^{2 \mu} \exp(- 2 \mu) T^{- 2 \mu }. $$
Hence, for all \ $\mu >0$, \ $w_\mu(T) \sim C_\mu T^{-2 \mu}$ \ holds.

On the other hand, we see that $f(t,\lambda)$ is monotonically 
decreasing. Hence, it takes the maximum value at $\lambda =0$:
$$ \sup_{\lambda >0} |f(t,\lambda)| \leq \T\frac{1}{2}\, T^2 . $$
Now, following the lines of the proof of \cite[Corollary~4.4]{EHN96}
(see also \cite[Remark~4.5]{EHN96}) we conclude that, with
$\frac{1}{T^2} = \alpha$, we have a regularization operator of
optimal order. \hfill \eop
\vskip0.5cm

If we compare the dynamic programming approach with
the Showalter me\-thod, they are quite similar with 
$T_{dyn}^2 \sim T_{Sw}$. Hence, to obtain the same 
order of convergence we only need $\sqrt{T_{Sw}}$
of the time for the Showalter method.

%---------------------------------------------------------------------------
\subsection{Regularization in the discrete case}

The dynamic programming principle allows us
to find an sequence of approximate solutions $\{u_k\}$ which
is a minimizer to a certain functional.

From regularization theory we are motivated to 
choose a functional which includes the norm of the residuals
$ \|Fu_k -y\|$. Since in general this will not necessarily 
yield a regularization, we include an additional term 
involving $u_{k+1}-u_{k}$.
Now analogous to the continuous case we want to minimize the functional
\begin{equation}\label{functional}
J(\{u_k\}_{k=0}^N) :=  \sum_{j=0}^{N} \|F u_k - y \|^2 +  \sum_{i=0}^{N-1} 
\|u_{k+1}-u_{k}\|^2 
\end{equation}
with respect to all sequences $\{u_k\}_{k=0}^N$ satisfying $u_0=0$.
The reason for choosing  the norm of the residuals is clear, since we want
to find an (approximate) solution to the equation $F u = y$.
The second term is important to obtain a regularization method, since
it controls the size of the steplength between two successive
iterations. 

At first sight it is not at all obvious that there is
a  constructive  method for minimizing
(\ref{functional}) with respect to all sequences $\{u_k\}_{k=0}^N$.
However, we show that the minimization problem
 can be treated within the framework of Subsection~\ref{discapp}.

Define $\epsilon_k$ as the $k$-th residual: 
\ $\epsilon_k := F u_k - y$, $k=0 \ldots, N$,
where $u_k$ is the solution we compute at the $k$-th iteration step.
The control is defined as $v_k = u_{k+1}-u_k$, $k = 0 \ldots N-1$.
As initial starting value we set $u_0 = 0$.
Hence we obtain the $k$-th iterate from the control variables by
\begin{equation}\label{iterate}
 u_k = \sum_{j=0}^{k-1} v_j. 
\end{equation}
From these definitions we obtain the following condition, 
which is trivially satisfied, when $v_k$ and $\epsilon_k$ are
defined in this way:
\begin{equation}\label{side}
 \epsilon_{k+1} = \epsilon_k + F v_k .
\end{equation}
Using the above notations, the minimization of (\ref{functional}), with 
initial condition $u_0 = 0$, is equivalent to the optimization problem 
in (\ref{eq:dcpr-hs}).

% \begin{equation}\label{functional1}
% J(z,u) = \sum_{j=0}^{N-1} (u_k,u_k) + 
% \sum_{j=0}^{N} (z_k,z_k)  
% \end{equation}
% subject to the conditions (\ref{side}) and the initial condition
% $$ z_0 = F x_0- y = - y $$

We now can use the results of Section~\ref{discapp} with $S = Q = R =
A = I$, $B = F$.
The dynamic programming principle yields  the iteration procedure
\begin{eqnarray}
 S_N & := & I \\
 K_k & := &  (F^*S_{k+1}F + I)^{-1} F^*S_{k+1}, \quad k = N-1 \ldots 0 \\
 S_k & := &  (I-FK_k)^*S_{k+1}(I-FK_k) + K_k^*K_k + I, \quad k = N-1,\ldots,0 
\end{eqnarray} 
If $K_k$, $S_k$ are computed, we obtain the control $v_k$ and the error 
$\epsilon_k$
from
\begin{eqnarray}
       \epsilon_0 & := & -y\;  \\
       v_k & = & - K_k \epsilon_k, \quad k=0,\ldots,N-1 \\
       \epsilon_{k+1} & = & \epsilon_k + F v_k = (I - F K_k) \epsilon_{k} . 
\end{eqnarray} 
The iterate $u_N$, which represents an approximation
to the solution, can be calculated from (\ref{iterate}).

Now we want to consider the mapping $y \to u_N$ as an iterative 
regularization operator where $N$ acts as regularization parameter.
% We will see, that this iteration can be treated within the standard 
% theory (as in \cite{EHN96}).
%
This mapping can be represented by filter functions $g_N$ using 
spectral theory, similar to the continuous case. The following
lemma serves as preparation for this purpose.
Let $ E_\lambda, F_\lambda$ be the spectral families of 
$F^*F$, $F F^*$.

\begin{lemma}
If $S_{k+1}$ has a representation as
$S_{k+1} = \int f_{k+1}(\lambda) d F_{\lambda}$,
with a continuous positive function $f_{k+1}$, then so has $S_{k}
 = \int f_{k}(\lambda) d F_{\lambda}$
and the following recursion formula holds:
\begin{equation}
f_k(\lambda) = 
\frac{ f_{k+1}(\lambda)( \lambda+1)  + 1} 
{ f_{k+1}(\lambda) \lambda + 1}  =
1+ \frac{ f_{k+1}} 
{ f_{k+1}(\lambda) \lambda + 1 } .
\end{equation}
\end{lemma} 

\noindent {\em Proof: }
We use the identity \ $ F^* f(FF^*) = f(F^* F) F^* $ \
\cite[formula (2.43)]{EHN96},
which holds for any piecewise continuous function $f$.
Since $f_{k+1}$ is positive, the inverse 
$\left(f_{k+1}(\lambda) \lambda + 1 \right)^{-1}$ exists, and
$$ K_k = 
\int \left(f_{k+1}(\lambda) \lambda + 1 \right)^{-1} 
f_{k+1}(\lambda) d E_{\lambda}\ F^* . $$
%
% Moreover,
% $$ K_k^* K_k = 
% F \int_{\sigma} \left(f_{k+1}(\lambda) \lambda + 1 \right)^{-2} 
% f_{k+1}(\lambda)^2 d E_{\lambda}  F^*  =
% \int_{\sigma} \left(f_{k+1}(\lambda) \lambda + 1 \right)^{-2} 
% f_{k+1}(\lambda)^2 \lambda d F_{\lambda}.    $$
% Similarly,
% $$ I- FK_k = \int_{\sigma} 1 -
%  \left(f_{k+1}(\lambda) \lambda + 1 \right)^{-1} 
% f_{k+1}(\lambda) \lambda d F_{\lambda} =
%  \int_{\sigma} 
%  \left(f_{k+1}(\lambda) \lambda + 1 \right)^{-1}  d F_{\lambda} $$
%
From the identity above and some basic algebraic manipulation we obtain
$$ S_k = 
 \int
 \left(f_{k+1}(\lambda) \lambda + 1 \right)^{-2} f_{k+1}
+ \left(f_{k+1}(\lambda) \lambda + 1 \right)^{-2} 
f_{k+1}(\lambda)^2 \lambda +1 
  d F_{\lambda} $$
%$$
%= \int_{\sigma} 
%\frac{f_{k+1}(\lambda) + \lambda f_{k+1}(\lambda)^2 +
%  (f_{k+1}(\lambda) \lambda + 1)^2} 
%{ \left(f_{k+1}(\lambda) \lambda + 1 \right)^{2}}
%  d F_{\lambda} $$
$$ = \int \frac{ (f_{k+1}(\lambda)( \lambda+1) + 1)} 
     { \left(f_{k+1}(\lambda) \lambda + 1 \right)} d F_{\lambda} = 
     \int 1 + \frac{f_{k+1}}{ f_{k+1}(\lambda) \lambda + 1 } d F_{\lambda} .
$$
\mbox{} \hfill \eop
\vskip0.5cm

By definition we have $S_N = I$,  $f_N$ obviously satisfies the hypothesis 
of the theorem with $f_N=1$ and hence, by induction, all $S_k$ have a 
representation via a spectral function $f_k$.

An obvious consequence of the recursion formula is the following recursion:
\begin{equation}\label{recursion}
h_{k}(\lambda) = 2 +\lambda - \frac{1}{h_{k+1}(\lambda)} ,
\end{equation}
with $h_{k}(\lambda) := \lambda f_{k}(\lambda) +1$ and the end condition
$h_{N}(\lambda) = \lambda + 1$.

% \begin{equation}\label{endcond}
% h_{N}(\lambda) = \lambda +1  
% \end{equation}

%\begin{equation}\label{recursion}
% h_{k} = 1+ \lambda \left(1+\frac{f_{k+1}}{h_{k+1}} \right)
% = 1 + \lambda + \frac{\lambda f_{k+1}}{h_{k+1}} 
% = 1+ \lambda + \frac{h_{k+1}-1}{h_{k+1}} = 
% 2+\lambda - \frac{1}{h_{k+1}} 
%\end{equation}

Now we want to find a filter function $g_N$ to express 
$ u_N = \int g_N(\lambda) d E_{\lambda} F^* y.$
Using the expression \ $I - FK_k = 
\int \left(f_{k+1}(\lambda) \lambda + 1 \right)^{-1}  d F_{\lambda}$
\ we conclude
$$ \epsilon_{k+1} = 
\int \left(f_{k+1}(\lambda) \lambda + 1 \right)^{-1}  d F_{\lambda} \epsilon_k =
\int \frac{1}{h_{k+1}} d F_{\lambda} \epsilon_k = 
-\int \frac{1}{\Pi_{i=1}^{k+1} h_i(\lambda)} d F_{\lambda} y $$
$$ v_k = 
-\int \frac{f_{k+1}(\lambda)}{h_{k+1}(\lambda)} d E_{\lambda} F^* \epsilon_k = 
% \int \frac{f_{k+1}(\lambda)}{h_{k+1}(\lambda)} 
% d E_{\lambda}
% \int \frac{1}{\Pi_{i=1}^{k} h_i(\lambda)} d E_{\lambda} F^*y =
\int \frac{f_{k+1}(\lambda)}
{ \Pi_{i=1}^{k+1} h_i(\lambda) } d E_{\lambda} F^* y $$
%
%\begin{equation}
% x_{k} = \sum_{i=0}^{k-1} 
%\int_{\sigma} \frac{f_{i+1}(\lambda)}{
% \Pi_{j=1}^{i+1} h_j(\lambda)} 
% d E_{\lambda} F^* y 
%\end{equation}
%
Now we replace $f_{i+1} = \frac{1}{\lambda} (h_{i+1}-1)$
and use (\ref{iterate}) to obtain
\begin{eqnarray}
 u_k & = & \sum_{i=0}^{k-1} 
\int_{\sigma} \frac{1}{\lambda} \frac{h_{i+1}(\lambda)-1}{
 \Pi_{j=1}^{i+1} h_j(\lambda)}  
 d E_{\lambda} F^* y \nonumber \\
 & = & \sum_{i=0}^{k-1} 
\int_{\sigma} \frac{1}{\lambda}  \left(
 \frac{1}{ \Pi_{j=1}^{i} h_j(\lambda)} -
   \frac{1}{ \Pi_{j=1}^{i+1} h_j(\lambda)} \right)
 d E_{\lambda} F^* y. \nonumber \\
& = & \int \frac{1}{\lambda} \left(
1-\frac{1}{\Pi_{j=1}^{k} h_j(\lambda)}  \right) d E_{\lambda} F^* y ,
\label{xrec}
\end{eqnarray}
where $h_k$ satisfies the backwards recursion formula (\ref{recursion})
and the end condition $h_{N}(\lambda) = \lambda +1$.

In particular, the $N$-th iterate, which is our approximate solution,
can be expressed as
$u_N = \int_{\sigma} g_N(\lambda) d E_{\lambda} F^* y$, with the
filter function
\begin{equation}\label{discfilter}
 g_N(\lambda) =  \frac{1}{\lambda} \left(
  1-\frac{1}{\Pi_{j=1}^{N} h_j(\lambda)}  \right) .
\end{equation}
\mbox{} \hfill \eop

The following theorem yields a representation for 
$g_N$ in Terms of Chebyshev polynomials.

\begin{theorem} 
Let $T_n(x)$ be the Chebyshev polynomial of the first kind of order $n$. Then
$$
g_N(\lambda) \ = \ \frac{1}{\lambda}
\left[ 1 - \left( \sqrt{\T\frac{\lambda}{4}+1} \right)
\left(T_{2N+1} \left(\sqrt{\T\frac{\lambda}{4}+1}\right) \right)^{-1}
\right] . $$
\end{theorem} 

\noindent {\em Proof:}
Define $p_i(\lambda) := \Pi_{k=N-i}^N h_k(\lambda)$, $i = 0\ldots N-1$. 
From the end condition for $h_N$ we find $p_0 = \lambda + 1$. Furthermore,
follows from (\ref{recursion})
\begin{equation}
\label{threerec} p_{i+1}(\lambda) = 
h_{N-i-1}(\lambda) p_i(\lambda) = 
(2+\lambda) p_i(\lambda) - 
\frac{p_i(\lambda)}{h_{N-i}(\lambda)} 
= (2+\lambda) p_i(\lambda) - p_{i-1}(\lambda),
\end{equation}
hence $p_i$ satisfies a three-term recursion.
From (\ref{recursion}) we see that $p_1 = \lambda^2 + 3 \lambda +1$.  
If we introduce $p_{-1}(\lambda) := 1$,
then the initial conditions $p_{-1}(\lambda)$, $p_{0}(\lambda)$
together with the three-term recursion (\ref{threerec}) 
completely determine $p_i$.

We prove the identity
$$ p_{N-1}(\lambda) \ = \ 
\frac{T_{2N+1} \left( \sqrt{\frac{\lambda}{4}+1} \right)}
{\sqrt{\frac{\lambda}{4}+1}} \ =: q_N(\lambda), \quad \forall N \geq 0 . $$
For $N = 0$ we have $p_{-1}(\lambda) = 1$
and, since $T_1(x) = x$, it follows $q_1 = 1$.  
Since $T_{3}(x) = 4 x^3 -3 x$ we find for $N=1$ that 
$q_2(\lambda) = \lambda +1 = p_1(\lambda)$.
Hence, the identity $p_{N-1}(\lambda)  = q_N(\lambda)$ holds
for $N=0,1$. Since  two initial conditions and the 
three-term recursion uniquely determine the
sequence $p_i(\lambda),q_i(\lambda)$ we only have to show
that  $q_i$  satisfies the same
recurrence relation as $p_i$.
Note that the following identity holds for
all $N \geq 1$ (cf. \cite[Page~132]{MOS66}):
$$ T_{2 N+3}(x) - T_{ 2N-1}(x) = 
2 T_{2N+1}(x) T_{2}(x) = 
2 T_{2N+1}(x) (2 x^2-1) . $$ 
Put $x = (\frac{\lambda}{4}+1)^{1/2}$ and multiply the identity by
$(\frac{\lambda}{4}+1)^{-1/2}$ we get
$$
\frac{T_{2 N+3}(\sqrt{\frac{\lambda}{4}+1})}
{\sqrt{\frac{\lambda}{4}+1}} -
\frac{T_{2 N-1}(\sqrt{\frac{\lambda}{4}+1})}
{\sqrt{\frac{\lambda}{4}+1}}
%
% & = &
% \frac{T_{2 N+1}(\sqrt{\frac{\lambda}{4}+1})}
% {\sqrt{\frac{\lambda}{4}+1}} 4 (\frac{\lambda}{4} +1) -2  \\
%
\ = \ \frac{T_{2 N+1}(\sqrt{\frac{\lambda}{4}+1})}
{\sqrt{\frac{\lambda}{4}+1}} (\lambda +2) .
$$
Thus $q_N$ satisfies $q_{N+1}(\lambda) = (\lambda +2) q_N - q_{N-1}$,
which is the same recurrence relation as $p_n$. Hence $q_{N} = p_{N-1}$.
\mbox{} \hfill \eop

\begin{corollary}
$g_N(\lambda)$ has the following representations:
\begin{eqnarray}
g_N(\lambda) & = & \frac{1}{\lambda}
\left( 1- 
\frac{\cosh\left( \mbox{\rm arcosh} (\sqrt{\frac{\lambda}{4}+1})\right)}
{\cosh\left( (2n+1)\mbox{\rm arcosh} \sqrt{\frac{\lambda}{4}+1}) \right)}
 \right) , \quad \lambda \geq 0 \label{cosh} \\
g_N(\lambda) & = & \frac{1}{\lambda}
\left( 1- \frac{1}{\sum_{m=0}^n \
\left( \begin{array}{c} 2 n +1 \\ 2 m \end{array} \right)
(\frac{\lambda}{4}+1)^{(n-m)} (\frac{\lambda}{4})^m} \right) . \label{poly}
\end{eqnarray}
\end{corollary}

\noindent {\em Proof:}
Equation (\ref{poly}) follows from the representation formula
for $T_{2n+1}$ (see \cite[Page~130]{MOS66}):
$$ T_{2n+1}(x) = 
\sum_{m=0}^n \left( \begin{array}{c} 2 n +1 \\ 2 m \end{array} 
\right) x^{2n+1-m} (x^2-1)^m . $$
For the identity (\ref{cosh}) we start with the
well-known representation (see \cite[Page~129]{MOS66})
$$ T_{n}(x) = \cos(n\ \mbox{arccos}(x)), \quad |x| \leq 1 $$
From $\cos(z) = \cosh(i z)$ and $\mbox{arcosh}(z) = i\ \mbox{arccos}(z)$ 
we get by analytic extension the identity
$$ T_n(x) \ = \ 
% \cos(n(\mbox{arccos}(z)) = 
% \lim_{z \in \C \to x } \cosh(n(\mbox{arcosh}(z)) = 
  \cosh(n\ \mbox{arcosh}(x)) , \ \ x \ge 1 . $$
From this representation (\ref{cosh}) follows, since $\lambda \geq 0$.
\mbox{} \hfill \eop
\vskip0.5cm

The next result concerns the regularization properties of the
proposed iterative method.

\begin{theorem}
The mapping $y \to u_N$ is a regularization operator, as $N \to \infty$.
\end{theorem} 

\noindent {\em Proof:}
We have to proof the similar properties for the filter function
$g_N(\lambda)$ as for the continuous case.

First of all, using L'H\^opital's rule we find
$$ \lim_{\lambda \to 0}
g_N(\lambda) = - \lim_{\lambda \to 0}
\frac{d}{d \lambda} \left(
\frac{\sqrt{\frac{\lambda}{4}+1}}
{T_{2 N+1} \sqrt{\frac{\lambda}{4}+1}} \right) = 
-\lim_{z \to 1} \frac{d}{d z}
\left( \frac{z} {T_{2 N+1}(z)} \right)
\left. \frac{1}{8 \sqrt{\frac{\lambda}{4}+1}} \right|_{\lambda = 0}
$$
$$
= \ -\frac{1}{8} \frac{T_{2 N+1}(1) - T'_{2N+1}(1)}{T_{2 N+1}(1)^2} \ = \
\frac{(2N+1)^2-1}{8}, $$
where we used $T_n(1) = 1$, $T_n'(1) = n^2$.
Hence $g_N(\lambda)$ can be extended continuously to $\lambda =0$,

The estimate $|\lambda g_N(\lambda) | \leq C$ 
reduces to  
$$   
\left| 1- 
\frac{\cosh\left( \mbox{\rm arcosh} (\sqrt{\frac{\lambda}{4}+1})\right)}
{\cosh\left( (2n+1)\mbox{\rm arcosh} \sqrt{\frac{\lambda}{4}+1}) \right)} 
\right| \ \leq \ C , $$
but, by the monotonicity of $\cosh$, it holds that
$0 \leq \frac{\cosh(x)}{\cosh((2n+1) x)} \leq 1$, as a consequence
the constant $C$ can be chosen $C=1$.

Finally,
$ \lim\limits_{N\to \infty } g_N(\lambda) \to \frac{1}{\lambda} $
holds, since $\lim\limits_{N\to \infty} \cosh((2N+1) x) = \infty$.
\hfill \eop
\vskip0.5cm

We now can proof the convergence rate result similar to the continuous case.
For this purpose, one has to estimate the obtained approximate solution for
the case of exact data and noisy data.

\begin{theorem}\label{th:convdisc}
Let $u_N$ be defined as above.
If the data are exact, $y = F u^\dagger$ and 
$u^\dagger$ satisfies a source condition (\ref{sourcecond}) for some
$\mu > 0$, then
\begin{equation} \label{rates2}
 \|u_N - u^\dagger\| \leq C_\mu N^{-2 \mu} .
\end{equation}

If the data are contaminated with noise,
$ \| y -y_\delta\| \leq \delta $ and $y = F u^\dagger$ 
with $u^\dagger$ satisfying (\ref{sourcecond}), then we have
constants $C_\mu,C$, independent of $N,\delta$, such that:
$$ \|u_{N,\delta} - u^\dagger\| \ \leq \ C_\mu N^{-2 \mu} + C \delta N. $$
The choice $N \sim \delta^{\frac{-1}{2 \mu +1}}$ 
yields the optimal order convergence rates 
\begin{equation} \label{rates3}
 \|u_{N,\delta} - u^\dagger\| \sim \delta^{\frac{2 }{2 \mu +1}} .
\end{equation}
\end{theorem}

\noindent {\em Proof:} 
We have to find an estimate for 
$$ |\lambda^\mu (1- \lambda g_N(\lambda)) | \leq w_\mu(N), \quad 
\forall \lambda \geq 0 . $$
Hence we need a bound for
$$ 
\xi(\lambda) \ := \ \frac{ \lambda^\mu 
\cosh\left( \mbox{\rm arcosh} (\sqrt{\frac{\lambda}{4}+1})\right)}
{\cosh\left( (2N+1)\mbox{\rm arcosh} \sqrt{\frac{\lambda}{4}+1}) \right)},
\quad \lambda \geq 0 . $$
We may transform the variables $x:= (\frac{\lambda}{4}+1)^{1/2}$,
$y = \mbox{\rm arcosh}(x)$ and, using $\cosh(y)^2 -1 = \sinh(x)^2$, we get
$$ \xi(\lambda(x(y))) \ = \
 \frac{ 4^\mu \sinh(y)^{2 \mu} \cosh(y) }
{\cosh\left( (2N+1) y) \right)} \ =: \ \zeta(y) , \quad  y \geq 0 . $$
For $y \geq 0$ we may use the addition theorems for $\cosh$:
$$
|\cosh((2N+1)y)| = |\cosh(2N y) \cosh(y) + \sinh( 2N y) \sinh(y)| $$
$$ =
|\cosh(y) \cosh(2 N y)| \left(1 + \tanh(2N y ) \tanh(y) \right)| \geq
|\cosh(y) \cosh(2 N y)|,
$$
and, with the estimate $\cosh(x) \geq \frac{1}{2}(\exp(x)+1)$, we get
$$ |\zeta(y)| \leq 
4^\mu \frac{ \sinh(y)^{2 \mu}}
{\cosh( 2 N y) }  
\leq  4^\mu 2 \frac{ \sinh(y)^{2 \mu}}
{\exp( 2 N y) +1 } :=  4^\mu 2 \ \eta(y) $$
Now differentiation yields the necessary condition for
a maximum of $\eta$: \ $\frac{ \mu}{ N}(1+ \exp(-x)) \ = \ \mbox{tanh}(x)$.
By  monotonicity we see that 
this equation has a unique solution $x^* > 0$ for $N > \mu$, which must be 
the maximum of $\eta(y)$, since $\eta(0) = 0$ and $\eta(\infty) = 0$.

Now express 
$ \sinh(x) = \frac{\tanh(x)}{\sqrt{1- \tanh(x)^2}} $,
use $\frac{1}{\exp(x)+1} \leq 1$, we get for $N > 2 \mu$
$$
\eta(x) \leq (\frac{\mu}{N})^{2 \mu} 
\frac{(1+ \exp(-x^*))^{\mu}} 
{\sqrt{1- \frac{ \mu^2}{ N^2}(1+ \exp(-x_*))^2} } \leq
C  \frac{1}{(2 N)^{2 \mu} } .
$$ 
Hence we get for all $\mu$ and $N > 2 \mu$
$$
\lambda^\mu|1- \lambda g_N(\lambda)| \leq 
C \frac{1}{(2 N)^{2 \mu}} ,
$$ 
which immediately yields (\ref{rates2}) (cf. \cite[Corollary~4.4]{EHN96}).

For a proof of (\ref{rates3}) we have to find an estimate
$$ g_N(\lambda) \leq C_N , \quad \forall  \lambda > 0 . $$
Using the same transformation as above, we have to bound
for all $y >0$,
$$ \phi(y):= 
\frac{\cosh((2 N+1)y) - \cosh(y)}{ \sinh(y)^2 \cosh((2 N+1) y)}
= \frac{2 \sinh((N+1)y) \sinh((N-1)y)}
{ \sinh(y)^2 \cosh((2 N+1) y)} $$
$$ \leq 2 
 \frac{\sinh(N y)^2}{\sinh(y)^2 \cosh(2N y)} 
\leq 2 \frac{\sinh((N+1)y)^2}{\sinh(y)^2( \cosh(N y)^2 + \sinh(N y)^2} $$
$$
\leq 
2 \left( \frac{\sinh((N+1)y)}{\sinh(y)( \cosh(N y))}\right)^2
=: 2 \psi(y)^2 .
$$
Now we may calculate  the derivative 
(using summation formula for $\sinh$, $\cosh$),
$$ 
\psi'(y) \ = \ 
\frac{N}{2} \left( \frac{\sinh(2 y)- \frac{1}{N} \sinh( 2 N y)}
{\sinh(y)^2\cosh(N y)^2} \right) . $$
Now by differentiation it is easy to see that
for positive $y$ the function  
$\sinh(2 y)- \frac{1}{N} \sinh( 2 N y)$ is 
 strictly monotonically decreasing  and it 
vanishes for $y=0$. Hence $\psi$ has negative derivative for $y >0$ 
and $\psi'(0) = 0$. Thus the maximum must be at $y = 0$.
By L'H\^opital's rule
$$ \psi(0) = \lim_{y \to 0} 
 \frac{\sinh((N+1)y)}{\sinh(y)} = N+1 . $$
Hence $|g_N(\lambda)| \leq 2 (N+1)^2 \leq C N^2$, with a constant $C$
independent of $N$. With the results of \cite[Theorem~4.3]{EHN96} the
proof is finished.
\mbox{} \hfill \eop

%---------------------------------------------------------------------------
% Section 4
%---------------------------------------------------------------------------
\section{Numerical experiments} \label{sec:num}

We are now concerned with the numerical realization
of the described algorithm. We consider the
discrete variant (\ref{eq:dcs},\ref{eq:dcu}) and
a discretization of the continuous algorithm
(\ref{eq:riccati},\ref{eq:cont-reg}).

The first one has a straightforward implementation.
For the continuous approach we use an explicit time-discretization
$Q'(t) \sim \tfrac{1}{\Delta t} (Q_{n+1}- Q_n)$.
Then Equation (\ref{eq:riccati}) becomes an
iterative procedure: (note that the Riccati-equation has
to be solved backwards in time)
\begin{eqnarray*}
Q_{n} & = & Q_{n+1} + \Delta t (I - Q_{n+1} F F^* Q_{n+1} ) ,
\quad n = N-1\ldots 0 \\
Q_N   & = & 0. 
\end{eqnarray*}

Equation (\ref{eq:cont-reg}) is discretized in a similar manner:
$$
u_{n+1} =  u_n  - \Delta t( F^* Q_n (F u_n - y)) ,
\quad n = 0 \ldots N-1
$$
together with some initial condition $u_0$.

A more efficient method is to use a recursion
for $B_n := F^*Q_n$. Since $Q_n$ is symmetric, then
\begin{equation}\label{eq:explicit_Ric}
 B_{n}  = B_{n+1} + \Delta t F^*(I - B^*  B_{n+1} ).
\end{equation}
Hence we get 
\begin{equation}\label{eq:explicit_it}
 u_{n+1} =  u_n  - \Delta t( B_n (F u_n - y)) .
\end{equation}

Since we used an explicit discretization scheme, the method
will be only stable if we bound the stepsize appropriately, e.g.,
$\Delta t \|F^*F\| \leq 1$.
The explicit discretization has the advantage that 
no matrix inversion is needed, by paying the price of 
a restricted stepsize.
A detailed analysis of the regularization properties of
this iterative scheme, in the spirit of Section~\ref{sec:regul},
is of course also possible.

As a benchmark problem we consider an integral equation
of the first kind:
$$ F u = \int_0^1 k(x,y) u(y) d y . $$ 
For a discretization of this operator, we split
the unit interval $I = [0,1]$ into $m$ subintervals and
discretize $u$ by using a uniform discretization with piecewise
linear, continuous splines on each subinterval
(also known as Courant-finite elements).
The integral is evaluated by the trapezoidal rule 
one each subinterval. As evaluation points for $x$ we used 
$x_i = i/m$, $i = 0, \ldots, m$.
This results in a $(m+1) \times (m+1)$ matrix equation:
\begin{equation} \label{eq:disceq}
 F_m u_m = y_m .
\end{equation}
We tested our algorithms with $F$ replaced by the discretized
version $F_m$.

We do not address the question how the discretization parameters
$m$ has to be related to the regularization parameter (the
iteration index in our case), but we simply consider 
the discretized equation as the given ill-posed problem.
Hence we use the Euclidean norm in $R^{m+1}$ on the discrete 
variables $u_m,y_m$. 

For our numerical test we used two different kernel functions $k(x,y)$:
\begin{eqnarray}
k_1(x,y) & := & \left\{ \begin{array}{cc}
                        (1-\frac{(x-y)^2}{0.1})^6 & \mbox{ if } 
(x-y)^2 \leq 0.1 \\
0 & \mbox{ else } \end{array} \right. \\
k_2(x,y) & := & \frac{1}{2 \sqrt{20}} \exp(-20 (x-y)^2). 
\end{eqnarray}
The first one is $6$-times continuously differentiable and 
hence leads to a mildly ill-posed problem. The second one
$k_2(x,y)$ is smooth, hence it leads to an exponentially
ill-posed problem.

We tested our methods for two exact solutions
$$ u_1^\dagger(x) \ := \ x(1-x) + \cos(20 x) , \quad u_2^\dagger(x) \ = \
\left\{ \begin{array}{cl} 1 & \mbox{if }\ 
0.3 \leq x \leq 0.5 \\
0 & \mbox{else} \end{array} \right. $$

We compared both algorithms with the Landweber-iteration and
the CG- method (see, e.g.,\cite{EHN96}). Throughout our numerical
experiments we used a discretization of $m = 300$.

\begin{figure}[ht]
\begin{center} 
\begin{minipage}{\textwidth}
\centerline{   \psfig{figure=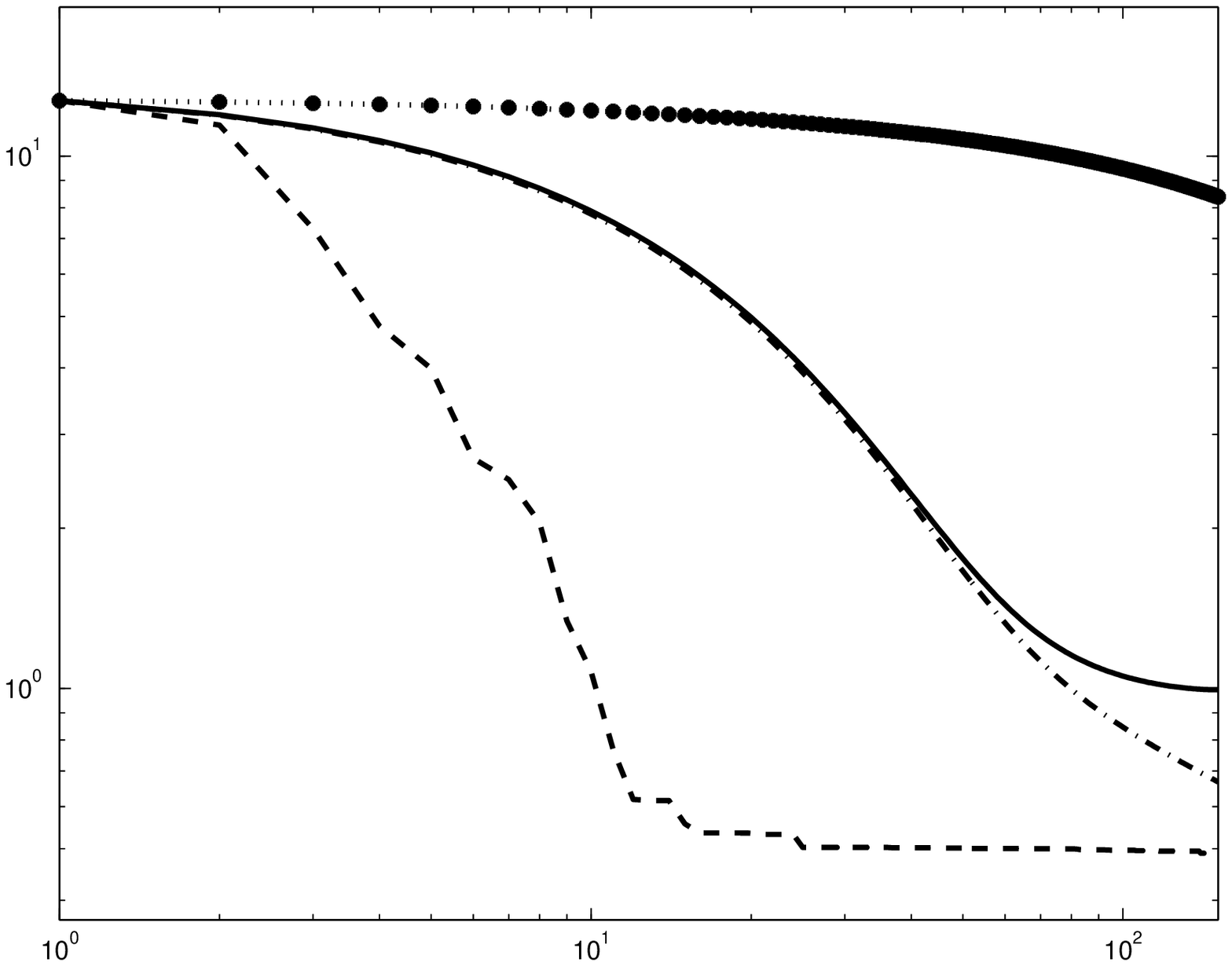,width=0.45\textwidth}
\hspace{0.5cm} \psfig{figure=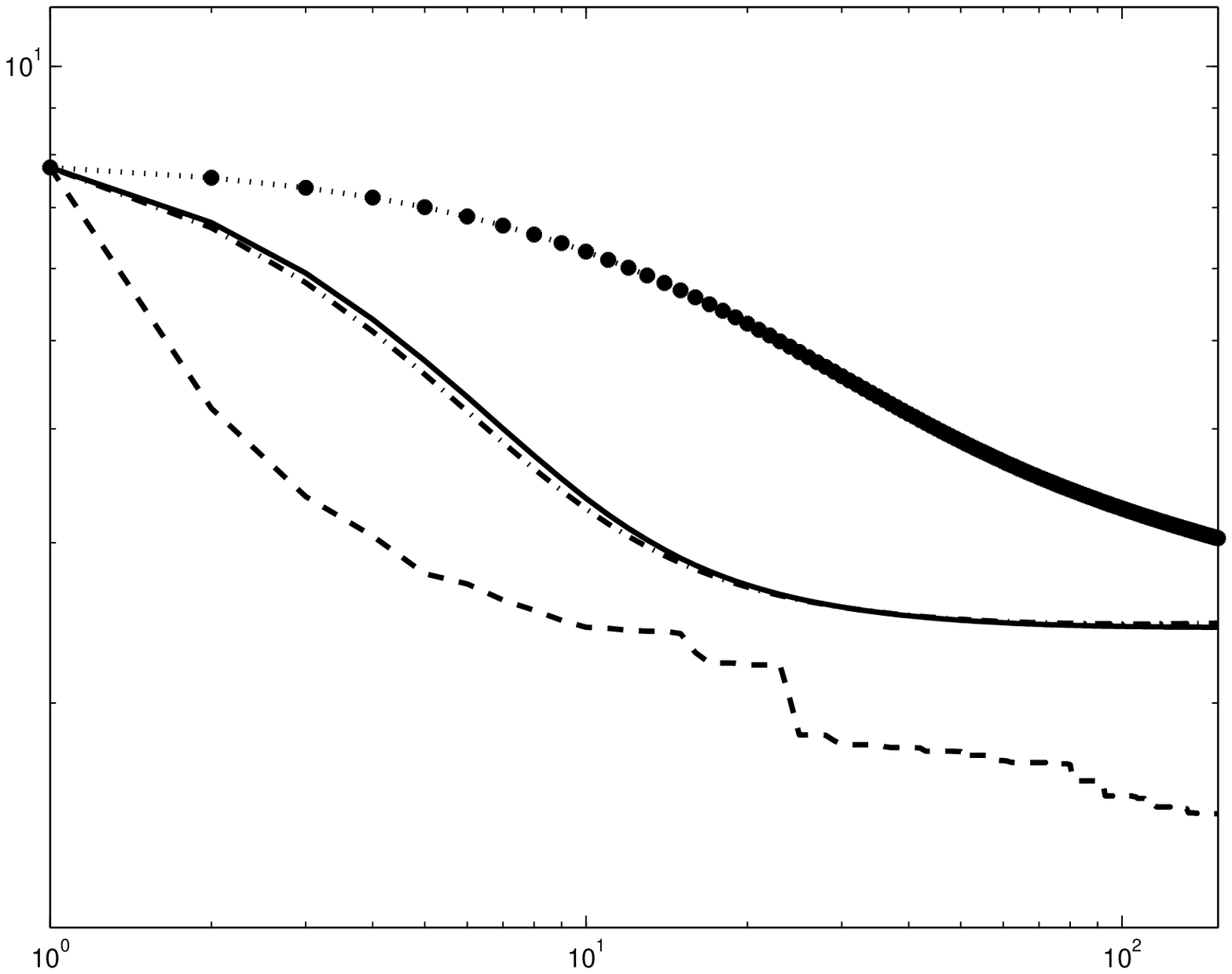,width=0.45\textwidth}}
\centerline{\hfil $k(x,y) = k_1(x,y)$, \ $u^\dagger=u_1^\dagger$
\hspace{2cm} $k(x,y) = k_1(x,y)$, \ $u^\dagger=u_2^\dagger$}
\vspace{0.5cm} 
\centerline{   \psfig{figure=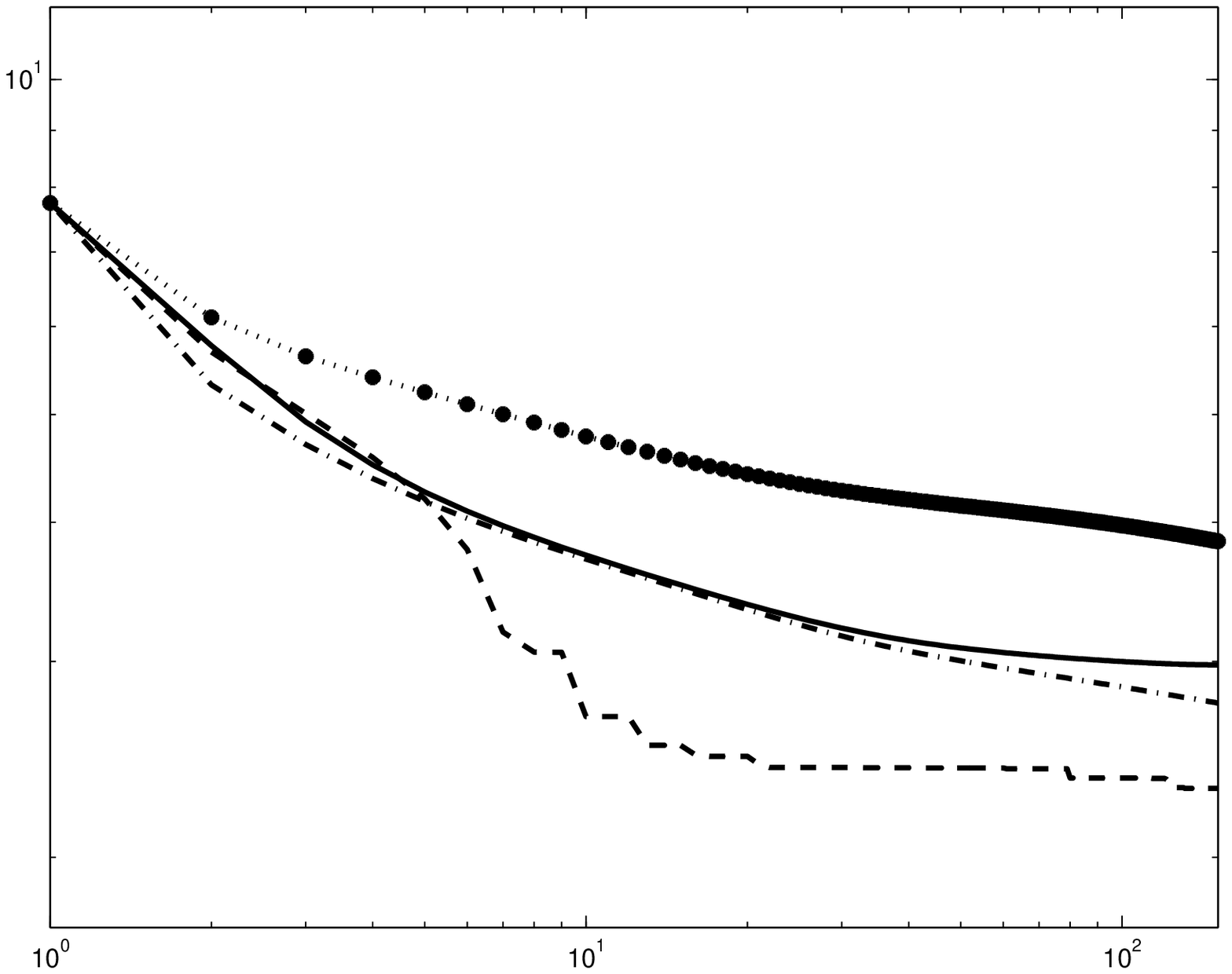,width=0.45\textwidth}%
\hspace{0.5cm} \psfig{figure=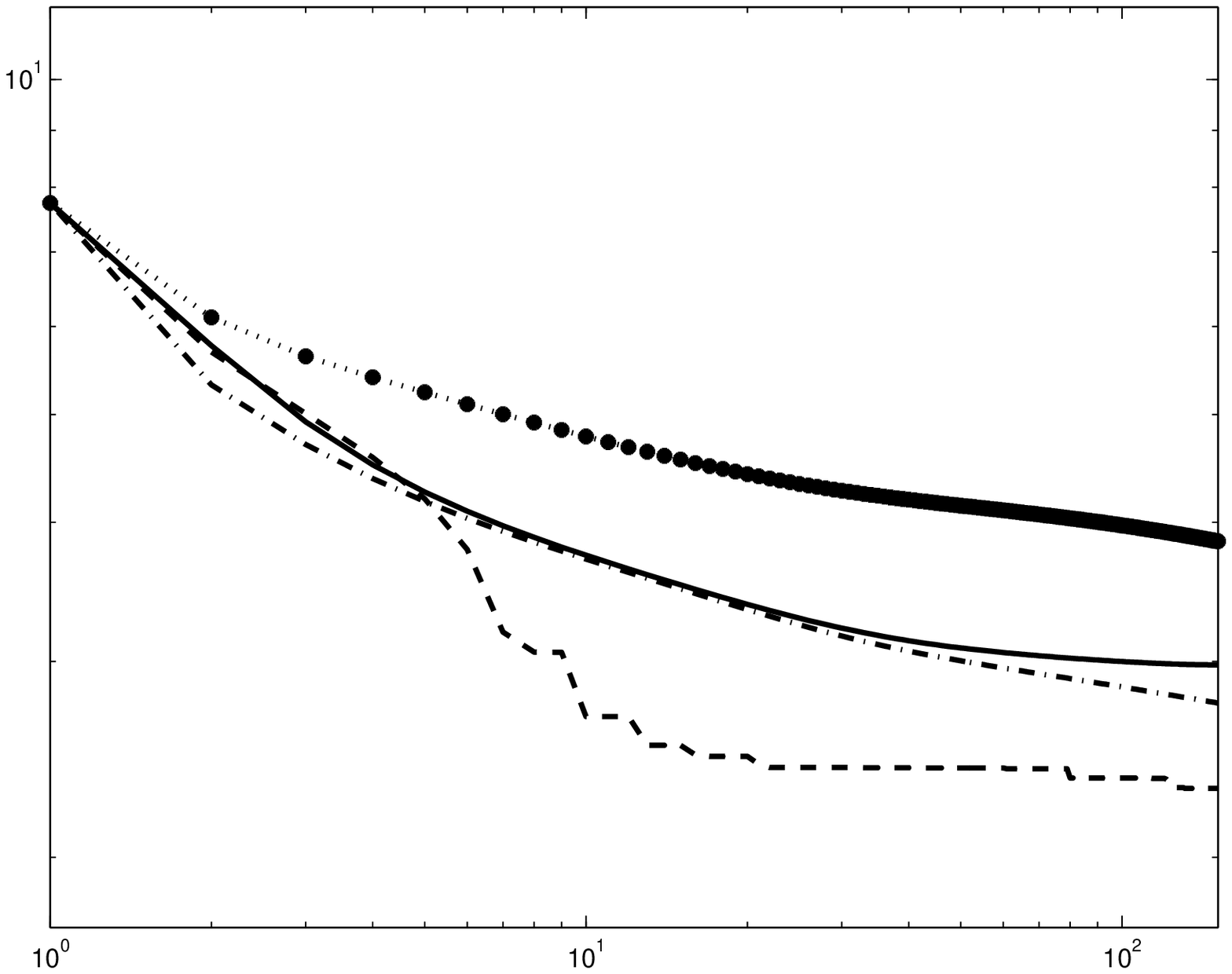,width=0.45\textwidth}}
\centerline{\hfil $k(x,y) = k_2(x,y)$, \ $u^\dagger=u_1^\dagger$
\hspace{2cm} $k(x,y) = k_2(x,y)$, \ $u^\dagger=u_2^\dagger$}
\end{minipage}
\end{center}\caption{\label{figure1} Evolution of the error 
$\|u_N - u^\dagger\|$ for exact data for all four algorithms.}
\end{figure}

Figure~\ref{figure1} shows the error $\|u_N-u^\dagger\|$ over the 
iteration index $N$ on a log-log scale for the four algorithms and 
the different choices of $u^\dagger$ and $k(x,y)$.
Here the full line corresponds to the discrete dynamic programming method, 
the dotted line to the Landweber iteration, the dashed-dotted to the 
continuous method with explicit time discretization, and the dashed 
line to the conjugate-gradient method.

\begin{figure}[ht]
\begin{center}
\begin{minipage}{\textwidth}
\centerline{   \psfig{figure=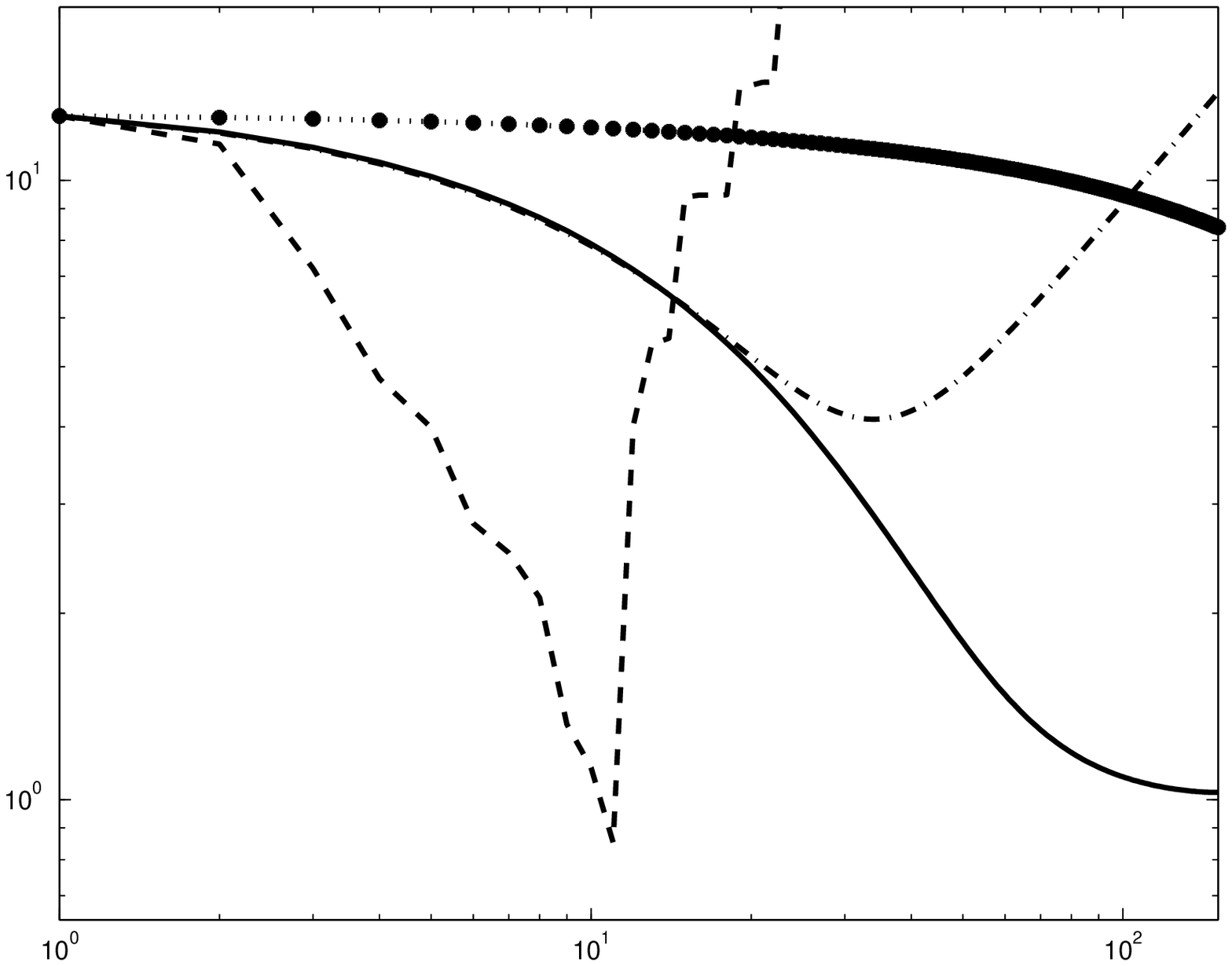,width=0.45\textwidth}%
\hspace{0.5cm} \psfig{figure=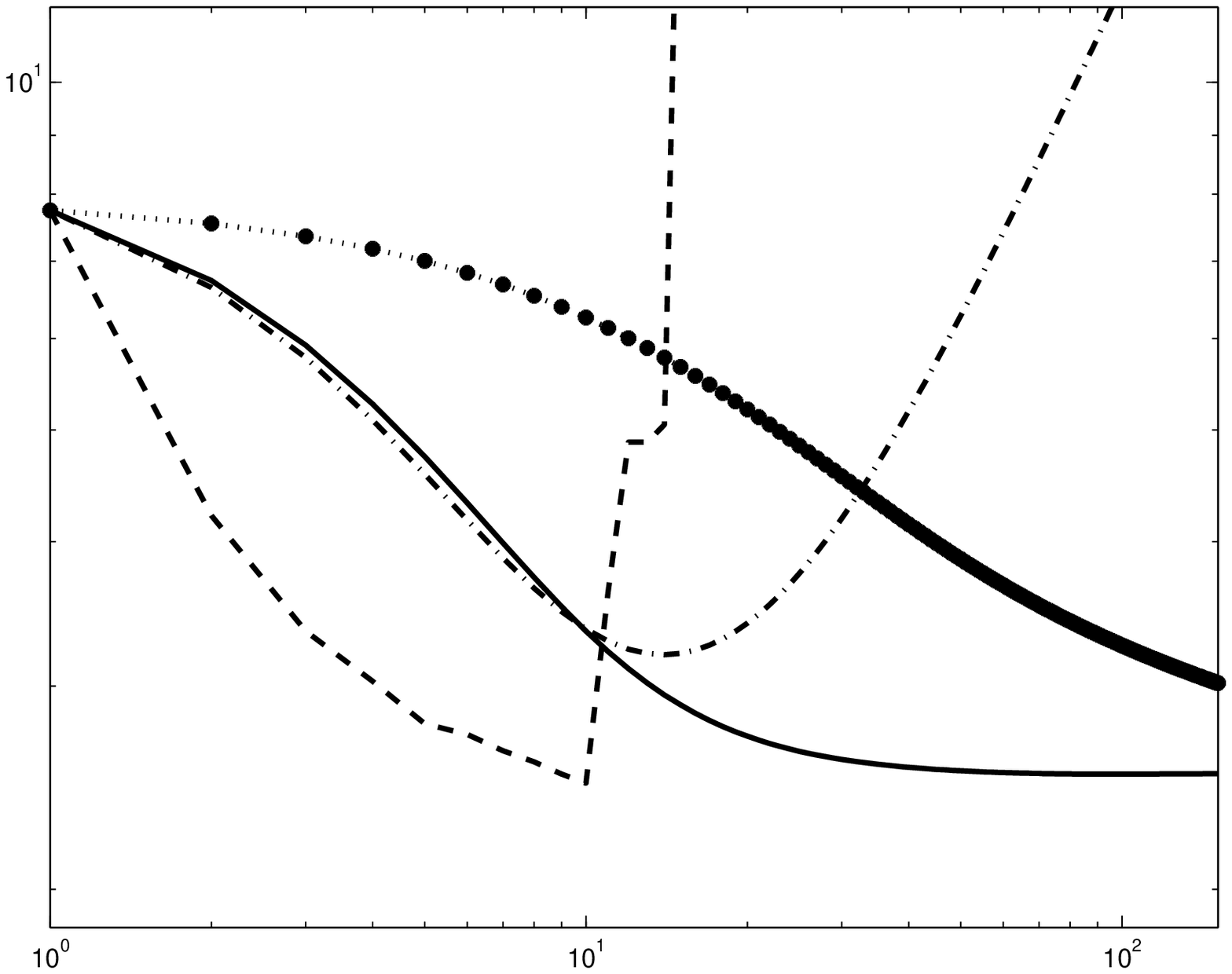,width=0.45\textwidth}}
\centerline{\hfil $k(x,y) = k_1(x,y)$, \ $u^\dagger=u_1^\dagger$
\hspace{2.0cm} $k(x,y) = k_1(x,y)$, \ $u^\dagger=u_2^\dagger$ \hfil}
\vspace{0.5cm} 
\centerline{   \psfig{figure=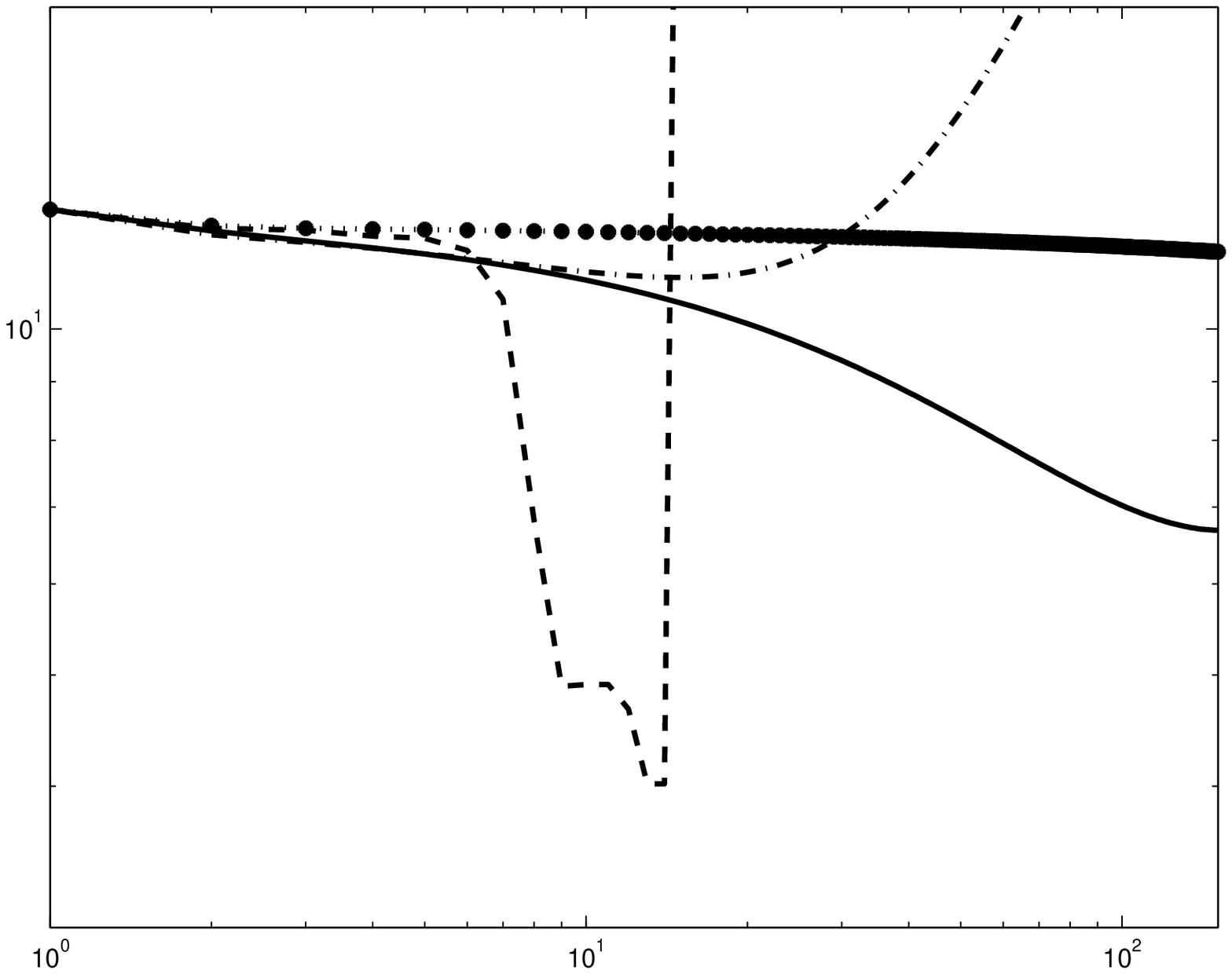,width=0.45\textwidth}%
\hspace{0.5cm} \psfig{figure=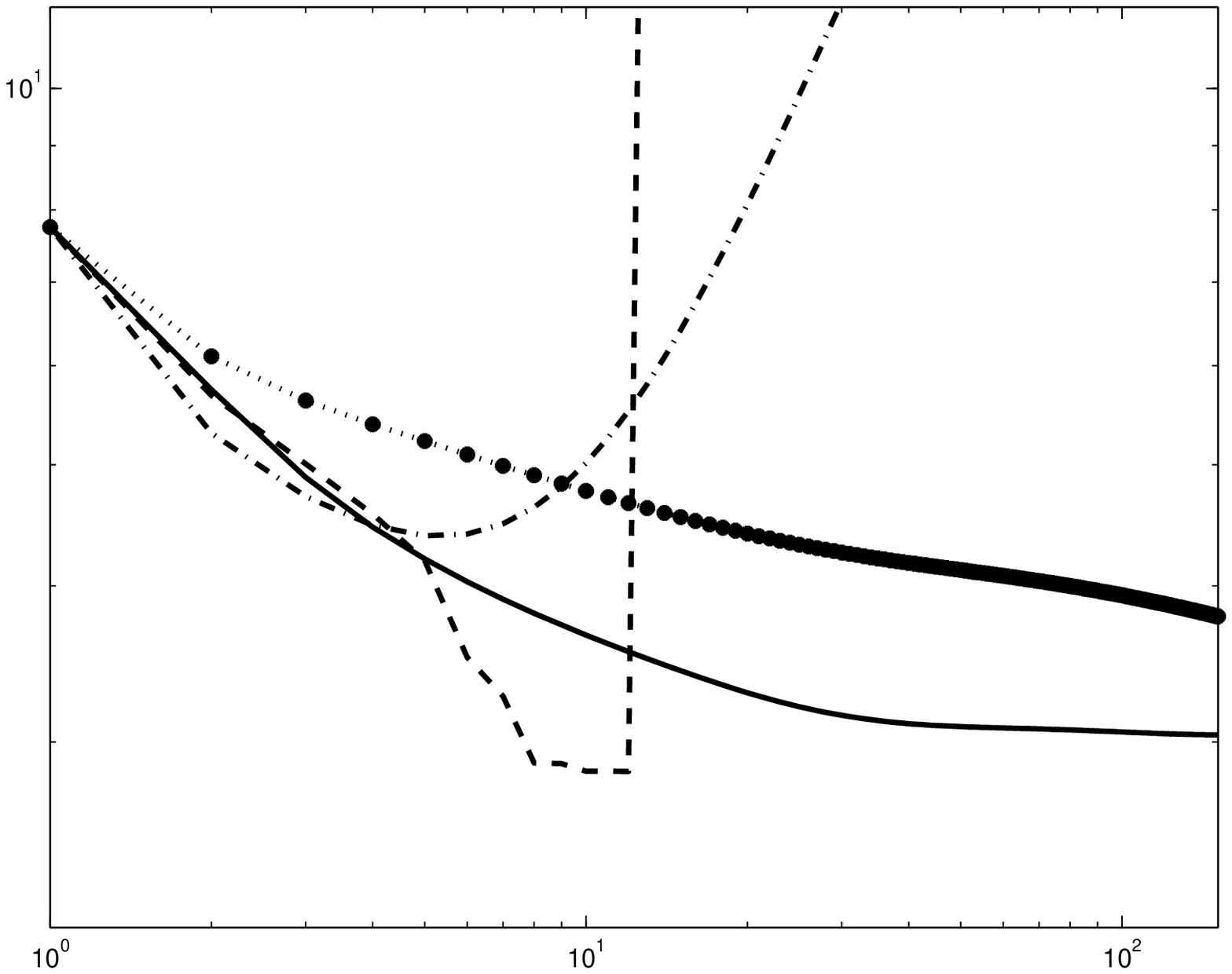,width=0.45\textwidth}}
\centerline{\hfil $k(x,y) = k_2(x,y)$, \ $u^\dagger=u_1^\dagger$
\hspace{2.0cm} $k(x,y) = k_2(x,y)$, \ $u^\dagger=u_2^\dagger$ \hfil}
\end{minipage}
\end{center}\caption{\label{figure3} Evolution of the error 
$\|u_N - u^\dagger\|$ for noisy data for all four algorithms.}
\end{figure}

Furthermore we contaminated the data with $10\%$ random noise.
The results are shown in Figure~\ref{figure3}.
Since in this case the iteration cannot converge,
a correct stopping criterion would be necessary. 
An a-priori stopping criterion was derived in 
Theorems~\ref{th:convcont} and~\ref{th:convdisc}.
Of course a-posteriori stopping criteria are more 
flexible. A more detailed analysis of these rules
(e.g., Morozov's discrepancy principle, or
the Engl-Gfrerer-type rules \cite{EHN96}) are 
out of the scope of this work.

We observe that the two methods based on dynamic programming techniques
are almost similar. Moreover these two methods have about 
the same convergence rates as the conjugate-gradient algorithm,
indicated by the same slope of the lines. 
This is confirmed by
theory, as the number of iterations $k$ to reach a certain noise level
$\delta$  under a source condition is $k \sim \delta^{\frac{-1}{2 \nu +1}}$ 
both for CG \cite[Thm~7.13]{EHN96}  and the dynamic 
programming techniques (Thm.~\ref{th:convdisc}), whereas for Landweber
iteration it is larger, namely $k \sim \delta^{\frac{-2}{2 \nu +1}}$.
Note also, that CG is only
a regularization method together with a discrepancy principle and is not
one in the sense of \cite{EHN96} if the noise level vanishes. Such a
phenomenon does not happen for the dynamic programming iterations.

Let us report on the overall costs of computation. Let $F$ be a matrix of
size $n \times m$. Then if $N$ time-steps (or iteration steps) are made,
the complexity for Landweber iteration and CG are ${\cal O}(n m N)$,
since only matrix-vector multiplications have to be performed.
The bottleneck for the dynamic programming iterations
\eqref{eq:dcs}-\eqref{eq:dcu}
and \eqref{eq:explicit_Ric}-\eqref{eq:explicit_it} is the Riccati equation.
Since in each step a matrix-matrix product has to be computed we end up
with an overall complexity for the implicit scheme
\eqref{eq:dcs}-\eqref{eq:dcu}
of ${\cal O}(n^2 m N + m^2 n N + n^3 N)$ and ${\cal O}(n^2 m N + m^2 n N)$
for the explicit one \eqref{eq:explicit_Ric}-\eqref{eq:explicit_it}.
This shows that these iterations have a complexity of at least one power 
higher than other iterations. If $n \sim m$, then the explicit and the
implicit dynamics iterations are even of comparable  complexity. In this case
the implicit version is to be favored as it has no stepsize restrictions. 

%---------------------------------------------------------------------------
%---------------------------------------------------------------------------
\section{Final remarks and conclusions}
In this article we combined control theory with abstract regularization
theory. We proposed iterative algorithms for solving linear inverse 
problems in Hilbert spaces and scrutinized their regularization properties. 
Our algorithms give rise to convergence and convergence rates under
the standard source conditions. The convergence properties are comparable
to a conjugate gradient method.

However, we have to admit, that in terms of computational complexity
our method is not really competitive with standard methods, as it 
involves matrix-matrix products in each iteration. 
On the other hand, the most costly part of our computation, the computation
of $Q(t)$ can be performed independent of the data. Hence, if for a 
fixed operator the same problem has to be solved with different data,
then $Q(t)$ only has to be computed once, e.g., by \eqref{eq:explicit_Ric} 
and the remaining iteration \eqref{eq:explicit_it}
involving the data is of similar complexity as the usual iteration methods. 
In this case our iterations are competitive with CG. 

Most of all we consider this work a good starting point into further
directions:
First of all it should be noticed that, if $Q_n$ is chosen constant, and
not computed by the Riccati equation, the continuous regularization 
method proposed in this paper reduces to a preconditioned Landweber 
iteration. Therefore, the dynamic programming regularization method can 
be considered as a generalization of the Landweber method. Since the
Landweber method is convergent we expect that solving the Riccati equation
is a numerical overkill. Instead one can think of solving the equation
inexact or using just a few number of steps of the Riccati iteration
to get a matrix $Q$, which can be used in a preconditioned Landweber (or CG)
iteration. 

Secondly, we expect that the real power of the combination of control theory
and regularization comes into play when considering dynamical inverse problems,
that is, if the data or the operator depend on time. In this case standard
iterations cannot be used, but the dynamic programming principle still can
be applied.

\section*{Acknowledgment}
The work of S.K. is supported by Austrian Science Foundation under
grant SFB F013/F1317; the work of A.L. is by CNPq, grant 306020/2006-8.

%---------------------------------------------------------------------------
%---------------------------------------------------------------------------

\end{document}